\documentclass[11pt]{article}

\usepackage{amsmath}
\usepackage{amssymb}
\usepackage{theorem}
\usepackage{euscript}
\usepackage{amscd}

\renewcommand{\title}[1]{
     \addvspace{3\baselineskip}  
     \begin{center} \LARGE \bf #1
     \end{center}
     \addvspace{2\baselineskip}}   

\renewcommand{\author}[1]{
     \addvspace{-1\baselineskip}  
     \begin{center} \large \sc #1
     \end{center}
     \addvspace{2\baselineskip}}   

\makeatletter

\def\section{%
        \@startsection{section}{1}{\z@}%
        {8ex plus 6ex minus 3ex}{\baselineskip}%
        {\normalfont\large\scshape\centering}%
        }

\renewcommand{\paragraph}[1]{{\par\removelastskip\vskip.5\baselineskip%
         \indent{\itshape{#1}}{\ifperiod.\else\global\periodtrue\fi}%
         \rm \ignorespaces}}

\let\goth=\mathfrak
\let\calligraphy=\mathcal

\def\NN{{\mathbb N}}

\def\PP{{\mathbb P}}
\def\QQ{{\mathbb Q}}

\def\ZZ{{\mathbb Z}}

\def\Jj{{\calligraphy J}}

\def\Oo{{\calligraphy O}}
\def\Pp{{\calligraphy P}}

\def\aaa{{\goth a}}
\def\BBB{{\goth B}}

\def\EEE{{\goth E}}

\def\ppp{{\goth p}}

\def\RRR{{\goth R}}

\def\VVV{{\goth V}}

\def\hbar{{\,\overline{\!h}}}

\def\qbar{{\,\overline{\!q}}}

\def\betabar{\operatorname{\overline{\beta}}}
\def\Ctilde{{\,\widetilde{\!C}}}

\def\tT{{\boldsymbol{T}}}

\def\Bl{\operatorname{Bl}}
\def\card{\operatorname{card}}

\def\ord{\operatorname{ord}}
\def\Pic{\operatorname{Pic}}

\let\lra=\longrightarrow

\def\ie{{\it i.e.}~}

\def\cf{{\it cf.}~}

\def\inv{^{-1}}

\let\phi=\varphi
\let\epsilon=\varepsilon

\newcommand{\floor}[1]{\left\lfloor#1\right\rfloor}
\newcommand{\ceil}[1]{\left\lceil#1\right\rceil}
\newcommand{\fpart}[1]{\left\langle#1\right\rangle}

\newcommand{\textq}[1]{\quad\text{#1}\quad}
\newcommand{\textql}[1]{\quad\text{#1}}

\def\theoname{Theorem}
\def\lemmaname{Lemma}
\def\propositionname{Proposition}
\def\notationname{Notation}
\def\corollaryname{Corollary}
\def\conjecturename{Conjecture}
\def\remarkname{Remark}
\def\remarksname{Remarks}
\def\examplename{Example}
\def\examplesname{Examples}
\def\definitionname{Definition}
\def\definitionsname{Definitions}
\def\notationname{Notation}

\def\proofname{Proof}

\def\Dquad{\hskip 0.6em plus .02em minus .2em}  
\def\Dpar{\belowdisplayskip=0pt\belowdisplayshortskip=0pt\par}

\def\bigpenalty{\interlinepenalty=\@M}
\def\smallpenalty{\interlinepenalty=100}

\newif\ifperiod \periodtrue 

\def\D@makemargins{%
  \labelsep=0pt
  \itemindent=0pt
  \labelwidth=0pt}

\def\D@restoremargins{%
  \labelsep=5pt
  \itemindent=0pt
  \leftmargin=5mm  
  \labelwidth=\leftmargin \advance\labelwidth by -\labelsep}

\def\th@Dindent{\hspace\parindent}
\def\th@Dheadingshape{\scshape}

\gdef\th@DthAndSuchtheo{%
  \D@makemargins%
  \def\@begintheorem##1##2{%
  \item[]\th@Dindent{\th@Dheadingshape ##1~\rm ##2.}\Dquad         
        \D@restoremargins}%
  \def\@opargbegintheorem##1##2##3{\def\next{##3}%
  \item[]\th@Dindent{\th@Dheadingshape ##1~\rm ##2\ifx\next\empty
  \else\ {\normalfont(##3)}\fi.}         
        \D@restoremargins}}

\gdef\th@DthAndSuchtheostar{%
  \D@makemargins%
  \def\@begintheorem##1##2{%
  \item[]\th@Dindent{\th@Dheadingshape ##1.}\Dquad     
        \D@restoremargins}%
  \def\@opargbegintheorem##1##2##3{\def\next{##3}%
  \item[]\th@Dindent{\th@Dheadingshape ##1\ifx\next\empty
  \else\ ##3\fi.}\Dquad         
        \D@restoremargins}}

\gdef\th@DthAndSuchliketheo{
  \D@makemargins%
  \def\@begintheorem##1##2{%
    \@latex@error{likethm: You must provide an argument in square brackets,
    though it may be empty [] !}%
    }%
  \def\@opargbegintheorem##1##2##3{%
        \def\next{##3}\ifx\next\empty\item[\th@Dindent]\else
        \item[]\th@Dindent{\th@Dheadingshape \next.}\Dquad\fi
        \D@restoremargins}}

\gdef\th@DdefAndSuch{%
  \D@makemargins%
  \def\@begintheorem##1##2{%
  \item[]\th@Dindent{\def@Dheadingshape ##1~\rm ##2.}\Dquad         
        \D@restoremargins}%
  \def\@opargbegintheorem##1##2##3{\def\next{##3}%
  \item[]\th@Dindent{\def@Dheadingshape ##1~\rm ##2\ifx\next\empty
  \else\ {\normalfont(##3)}\fi.}         
        \D@restoremargins}}

\gdef\th@DdefAndSuchStar{%
  \D@makemargins%
  \def\@begintheorem##1##2{%
  \item[]\th@Dindent{\def@Dheadingshape ##1.}\Dquad     
        \D@restoremargins}%
  \def\@opargbegintheorem##1##2##3{\def\next{##3}%
  \item[]\th@Dindent{\th@Dheadingshape ##1\ifx\next\empty
  \else\ ##3\fi.}\Dquad         
        \D@restoremargins}}

\def\th@Dheadingshape{\scshape}
\def\def@Dheadingshape{\itshape}

\theoremstyle{DthAndSuchliketheo}
\theorembodyfont{\bigpenalty\itshape}   
\newtheorem{likethm}{}
\theorembodyfont{\rmfamily}  

\theoremstyle{DthAndSuchtheostar}
\theorembodyfont{\bigpenalty\itshape}
\newtheorem{thm*}{\theoname}
\newtheorem{lem*}{\lemmaname}
\newtheorem{pro*}{\propositionname}
\newtheorem{cor*}{\corollaryname}
\newtheorem{conjecture*}{\conjecturename}
\theorembodyfont{\smallpenalty\rmfamily}
\newtheorem{notation*}{\notationname}
\newtheorem{exa*}{\examplename}
\newtheorem{examples*}{\examplesname}
\theoremstyle{DdefAndSuchStar}
\theorembodyfont{\smallpenalty\rmfamily}
\newtheorem{definition*}{\definitionname}
\newtheorem{definitions*}{\definitionsname}
\newtheorem{rem*}{\remarkname}
\newtheorem{remarks*}{\remarksname}
\theoremstyle{DthAndSuchtheo}
\theorembodyfont{\bigpenalty\itshape}
\newtheorem{thm}{\theoname}[section]
\newtheorem{lem}[thm]{\lemmaname}
\newtheorem{pro}[thm]{\propositionname}
\newtheorem{cor}[thm]{\corollaryname}

\theorembodyfont{\smallpenalty\rmfamily}
\newtheorem{exa}[thm]{\examplename}

\theoremstyle{DdefAndSuch}
\theorembodyfont{\smallpenalty\rmfamily}
\newtheorem{definition}[thm]{\definitionname}
\newtheorem{rem}[thm]{\remarkname}

\newtheorem{remarks}[thm]{\remarksname}

\theoremstyle{DthAndSuchtheo}               
\theorembodyfont{\itshape}

\newcommand{\proof}[1][]{{\par\removelastskip\vskip.6\baselineskip   
    \noindent\th@Dindent\def\next{#1}%
    {\itshape\proofname\ifx\next\empty\else\next\fi\ifperiod.%
      \else\global\periodtrue\fi\Dquad}%
    \clubpenalty=5000\rm\ignorespaces}\setcounter{step}{0}}

\newcounter{step}

\newcommand{\likeproof}[1][]{{\par\removelastskip\vskip.6\baselineskip
    \noindent\th@Dindent\def\next{#1}%
    {\itshape\ifx\next\empty\else\next\fi\ifperiod.%
      \else\global\periodtrue\fi\Dquad}%
    \clubpenalty=5000\rm\ignorespaces}\hspace{-2pt}\setcounter{step}{0}}

\def\qed{{\ifmmode\hskip 6mm plus 1mm minus 3mm{$\square$}
    \else
    \hfil\penalty50\hskip1em\null\nobreak\hfil
    {\hfill $\square$\parfillskip=0pt\finalhyphendemerits=0
      \let\par=\endgraf\par}
    \fi
    \Dpar\penalty-150\vskip.6\normalbaselineskip}}

\makeatother

\usepackage{pstricks}
\usepackage{pst-node,pst-tree}
\usepackage[arrow,curve,matrix,tips,frame]{xy}

\begin{document}

\title{Jumping numbers of a unibranch curve on a smooth surface}
\author{Daniel Naie}

\begin{abstract}
A formula for the jumping numbers of a curve unibranch at a singular
point is established.  The jumping numbers are expressed in terms of
the Enriques diagram of the log resolution of the singularity, or
equivalently in terms of the canonical set of generators of the
semigroup of the curve at the singular point.
\end{abstract}

The jumping numbers of a curve on a smooth complex surface are a
sequence of positive rational numbers revealing information about the
singularities of the curve.  They extend in a natural way the
information given by the log-canonical threshold, the smallest jumping
number (see \cite{EiLa} for example).  They are periodic, completely
determined by the jumping numbers less than $1$, but otherwise
difficult to compute in general, even if a set of candidates is easy
to provide, \cf \cite[Lemma~9.3.16]{La}.

The aim of this paper is to give a formula for the jumping numbers of
a curve {\it unibranch} at a singular point.  A curve $C$ will be said
to be unibranch at a point $P$, if the analytic germ of $C$ at $P$ is
irreducible.  The formula is expressed in terms of the Enriques
diagram associated to the singularity, or equivalently (see 
Theorem~\ref{th:main}) in terms of a minimal set of generators
$(\betabar_0,\betabar_1,\ldots,\betabar_g)$ of the
semigroup $S(C,P)$ of $C$ at $P$:
\begin{equation} \label{eq:main}
  \{\text{jumping numbers}<1\} = 
  \bigcup_{j=1}^g \frac{1}{[m_j,\betabar_j]}\, 
  R^{m_{j+1}}\bigg(\frac{m_j}{m_{j+1}},\frac{\betabar_j}{m_{j+1}}\bigg),
\end{equation}
where the $m_j$ are defined below.  Here
\[
  R^m(p,q) = \bigcup_{k=0}^{m-1}(kpq+R(p,q))
\]
and
\[
  R(p,q) = R^1(p,q) = \{ap+bq \mid a,b\in\NN^\ast,\, ap+bq<pq\}.
\]
The semigroup is defined by 
$S(C)=\{\ord_P s\mid s\in\Oo_{C,P}\}$, the order of the local
section $s$ being computed using a normalization of $C$.  It is finitely
generated and a minimal set of generators
$(\betabar_0,\betabar_1,\ldots,\betabar_g)$ is constructed as follows
(see \cite[Theorem~4.3.5]{Wa}): $\betabar_0$ is the least element of
$S(C)$; set $m_1=\betabar_0$; $\betabar_j$ is the least element of
$S(C)$ not divisible by $m_j$ and $m_{j+1}=\gcd(m_j,\betabar_j)$.

To prove (\ref{eq:main}) we use the notion of {\it~relevant divisors}
of the minimal log resolution of $C$ at $P$, notion introduced in 
\cite{SmTh}, and previously in \cite{FaJo} from the point of view of
valuations corresponding to Puiseux exponents: a relevant
divisor is an irreducible exceptional divisor that intersects at least
three other components of the total transform of $C$ through the
resolution.  When $C$ is unibranch at $P$, we show that the relevant
divisors account for all the jumping numbers.  This is the content of
Proposition~\ref{p:unibranch} and represents the key step of the
proof.  In general, if the curve is not unibranch, there are jumping
numbers that are not contributed by any relevant divisor (see
\cite[Example~2.2]{SmTh}).  Using Proposition~\ref{p:unibranch} and
the Enriques diagram associated to the minimal log resolution of $C$
at $P$, we compute the jumping numbers contributed by the relevant
divisors, and hence all the jumping numbers in
Theorem~\ref{th:mainEnriques}.  Finally, Theorem \ref{th:main} follows
as a consequence of Theorem \ref{th:mainEnriques} and of the
equivalence between the Enriques diagram and the semigroup $S(C)$.

The construction of the Enriques diagram as well as the definition of
the jumping numbers are recalled in \S\,\ref{s:one}.  The proof of
Theorem~\ref{th:mainEnriques} is given in \S\,\ref{s:three} together
with some necessary technical lemmas, whereas
Proposition~\ref{p:unibranch} is established in the last section.  The
explicit equivalence between the Puiseux characteristic, and hence the
semigroup, and the Enriques diagram of a unibranch curve is presented
in Theorem~\ref{th:enriques}.

In \cite{Ja}, Tarmo J\"arviletho obtained recently an explicit
description of the jumping numbers of a simple complete ideal $\ppp$
in a two dimensional regular local ring.  The jumping numbers are
expressed in terms of the Zariski exponents of the ideal.  Moreover
(see \cite[Proposition~9.2.8]{La} and also \cite[Theorem~9.4]{Ja}) the
jumping numbers $<1$ of the ideal $\ppp$ coincide to those of the
unibranch plane curve corresponding to a general element of $\ppp$,
and they amount to the jumping numbers given in (\ref{eq:main}).

If the unibranch curve is characterized by a single Puiseux exponent
$q/p$, with $\gcd(p,q)=1$, or equivalently if the semigroup is
generated by $p$ and $q$, the jumping numbers 
\[
  \frac{ap+bq}{pq}<1,\quad a,b\in\NN^\ast,
\]
were computed by L.~Ein in \cite{Ei}, or by J.~A.~Howald in \cite{Ho}
as a particular case of his formula for the multiplier ideals of
monomial ideals.

It is a pleasure to acknowledge James Alexander's criticism on
questions of presentation and proportion after having read a
preliminary version of the paper, and also the friendly and useful
talks I have with him, with Michel Granger and Adam Parusinski.

\section{Preliminaries and notation}
\label{s:one}

In this section we recall the definition of the jumping numbers and
introduce the Enriques diagram associated to a minimal log resolution
of a curve at a singular point.  The diagram will be used to perform
the calculations of the jumping numbers.

\subsection{Log resolutions and Enriques diagrams}
\label{ss:enriquesDiagrams}
Let $C$ be a curve on a smooth surface with an isolated singularity at
$P$.  A minimal log resolution of $C$ at $P$ is the composition
$\mu:Y\to X$ of blowings up such that $\mu$ gives an isomorphism
$Y\smallsetminus\mu\inv(P)\to X\smallsetminus\{P\}$, the strict
transform $\Ctilde$ of $C$ is smooth, the support of the total
transform $\mu^\ast C$ has normal crossings and the number of blowings
up is minimal with these properties.

Let $Y=Y_{r+1}\to Y_r\to\cdots\to Y_1=X$ be a decomposition of
$\mu:Y\to X$ into successive blowings up with
$Y_{\alpha+1}=\Bl_{P_\alpha}Y_\alpha$.  Each point $P_\alpha$ is
infinitely near to $P=P_1$ and has an associated exceptional divisor%
\footnote{The terms exceptional divisor and exceptional curve
will be used indifferently in the sequel.}
$E_\alpha$ on $Y_{\alpha+1}$.  Its strict
transform on $Y$ will also be denoted by $E_\alpha$ and its total
transform on $Y$ will be denoted by $W_\alpha$.
The strict transforms $E_\alpha$ and the the total transforms
$W_\alpha$ form two bases of the $\ZZ$-module
$\Lambda_C=\bigoplus_\alpha\ZZ E_\alpha\subset\Pic Y$.  In particular
\[
  \mu^\ast C = \Ctilde + D 
  = \Ctilde + \sum_\alpha e_\alpha E_\alpha
  = \Ctilde + \sum_\alpha w_\alpha W_\alpha.
\]

\begin{rem}
  \label{r:theWeights}
The weight $w_\alpha$, the coefficient of the total transform
$W_\alpha$ in the divisor $\mu^\ast C-\Ctilde$, is the multiplicity of
the corresponding strict transform of $C$ at $P_\alpha$.
\end{rem}

If the curve $C$ is unibranch, then for any $\alpha$, the strict
transform of $C$ on $Y_\alpha$ has a unique singular point $P_\alpha$.
In particular there is a unique log resolution of the singular point
$C$.  In the general case the resolution is not unique; the ordering
of the exceptional divisors $E_\alpha$, or equivalently of the points
$P_\alpha$ might vary.  Nevertheless, the ordering of the points is
compatible with the partial order of the infinitely near points.  If
$\alpha<\beta$, then either $P_\beta$ is infinitely near to $P_\alpha$
or there is $\gamma<\alpha$ such that $P_\alpha$ and $P_\beta$ are
infinitely near to $P_\gamma$.

The combinatorics of the configuration of the exceptional curves
$E_\alpha$ on $Y$, or equivalently the geometric relation between the
infinitely near points $P_\alpha$, is encoded in the notion of {\it
proximity}: a point $P_\beta$ is said to be {\it proximate} to
$P_\alpha$ if $P_\beta$ lies on the strict transform of $E_\alpha$ on
$Y_\beta$.  A point that is infinitely near is always proximate to at
most two other points.  A point is said to be {\it free} if it is
proximate to exactly one other point and {\it satellite} if it is
proximate to two infinitely near points.

\begin{rem}
  \label{r:theOrder}
If $P_\beta$ is proximate to both $P_\alpha$ and $P_{\alpha'}$ and if
$\alpha<\alpha'$, then $P_{\alpha'}$ is infinitely near to $P_\alpha$.
Moreover, if the curve is unibranch, then $P_\beta$ is always
proximate to $P_{\beta-1}$.
\end{rem}

A convenient way to present the proximity relations is given by the
{\it proximity matrix} $\Pi=||p_{\alpha\beta}||$, where
$p_{\alpha\alpha}=1$ for any $\alpha$ and $p_{\alpha\beta}$ equals
$-1$ if $P_\beta$ is proximate to $P_\alpha$ and $0$ if not.  The
proximity matrix is upper unitriangular by the previous remark and
represents at the same time the decomposition matrix of the strict
transforms in terms of the total transforms on $Y$.  A simple but
useful consequence of this remark is the following lemma.

\begin{lem}
  \label{l:star}
Let $D=\sum_\alpha e_\alpha E_\alpha=\sum_\alpha w_\alpha W_\alpha$ be
the divisor associated to a log resolution of $C$ at $P$. 
If $P_\beta$ is a satellite point
proximate to $P_{\alpha'}$ and $P_{\alpha''}$, then
\[
  e_\beta = e_{\alpha'}+e_{\alpha''}+w_\beta.
\]  
\end{lem}

\proof
Use the relation 
$E_\alpha=\sum_\beta p_{\alpha\beta}W_\beta=
W_\alpha-\sum_{P_\beta \text{ proximate to }P_\alpha}W_\beta$
to express the coefficients $w_\beta$ in terms of the coefficients
$e_\alpha$. 
\qed

The resolution data of a curve $C$ at a singular point $P$ has been
encoded by Enriques in an appropriate weighted tree diagram now called
the {\it Enriques diagram} (see \cite{EnCh, SeKn, Ca, Ev}).  The tree
graphically represents the proximity relations of the infinitely near
points.

\begin{definition*}
An {\it Enriques tree} is a couple $(T,\epsilon_T)$, where
$T=T(\VVV,\EEE)$ is an oriented tree (a graph without loops) with a
single {\it root}, with $\VVV$ the set of vertices and $\EEE$ the set
of edges, and where $\epsilon_T$ is a map
\[
  \epsilon_T:
  \EEE \to \{\text{`slant'}, \text{`horizontal'},\text{`vertical'}\}.
\] 
fixing the graphical representation of the edges.
An {\it Enriques diagram} is a weighted or labeled Enriques tree. 
\end{definition*}

\begin{definition*}
Let $T$ be an Enriques tree.  A horizontal (respectively vertical)
{\it $L$-shape branch} of $T$ is a path of length $\geq1$ such that
all edges but the first, are horizontal (respectively vertical)
through $\epsilon_T$.  An $L$-shape branch is proper if it contains at
least two edges.  A {\it maximal $L$-shape branch} is an $L$-shape
branch that cannot be continued to a longer one.  
\end{definition*}

The construction of the Enriques tree associated to a log resolution
of $C$ at $P$ is as follows.  The set of vertices is
$\VVV=\{P_1,\ldots,P_r\}$, \ie the set of infinitely near points; the
root of the tree is the proper point $P$.  There is an edge starting
at $P_\alpha$ and ending at $P_\beta$ if and only if $P_\beta$ is
proximate to $P_\alpha$ and, either $P_\beta$ is free, or $P_\beta$ is
satellite, proximate to $P_\alpha$ and $P_{\alpha'}$ and
$\alpha>\alpha'$.  There is an $L$-shape branch that starts at
$P_\alpha$ and ends at $P_\beta$ if and only if $P_\beta$ is proximate
to $P_\alpha$; there is either a horizontal or a vertical edge that
ends at $P_\alpha$ if and only if $P_\alpha$ is satellite.  To
normalize the shape of the tree it is assumed that an edge that starts
at a free point and ends at a satellite point is horizontal.  The
weights $w_\alpha$ are given by the coefficients of the total
transforms in $\mu^\ast C=\Ctilde+\sum_\alpha w_\alpha W_\alpha$.

It is to be noticed that $E_\alpha$ and $E_\beta$ intersect on $Y$ if
and only if there is a maximal $L$-shape branch that starts at
$P_\alpha$ and ends at $P_\beta$.

\begin{exa}[Definition of $T_{p,q}$]
  \label{ex:tpq}
Let $p<q$ be relatively prime positive integers and let $C$ be defined
locally at $P$ by $x^p-y^q=0$.  The Enriques tree $T_{p,q}$ associated
to the minimal log resolution of $C$ at $P$ is defined as follows.
Consider the Euclidean algorithm: $r_0=a_1r_1+r_2$, \ldots,
$r_{m-2}=a_{m-1}r_{m-1}+r_{m}$ and $r_{m-1}=a_mr_m$, with $r_0=q$ and
$r_1=p$.  Set
\[
  \VVV=\{P_\alpha \mid 1\leq\alpha\leq a_1+\cdots+a_m=r\} 
\]
and
\[
  \EEE=\{[P_\alpha P_{\alpha+1}] \mid
  1\leq\alpha\leq a_1+\cdots+a_m-1\}.
\]  
The map $\epsilon$ is locally constant on the $a_j$ edges 
$[P_\alpha P_{\alpha+1}]$ with 
$a_1+\cdots+a_{j-1}+1\leq\alpha\leq a_1+\cdots+a_j$.  The first
constant value of $\epsilon$---on the first $a_1$ edges---is `slant'.
The other constant values are alternatively either `horizontal' or
`vertical', starting with `horizontal'.  The Enriques trees $T_{5,7}$
is
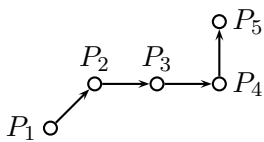
\begin{center}
\begin{pspicture}(0,0)(10,2.25)
  \psset{arrows=->,radius=.1,unit=5ex}
  \rput(5,1){
    \Cnode(-.71,-.71){1}
    \Cnode(0,0){2}
    \Cnode(1,0){3}
    \Cnode(2,0){4}
    \Cnode(2,1){5}
    \ncline{1}{2}
    \ncline{2}{3}
    \ncline{3}{4}
    \ncline{4}{5}
    \uput{\labelsep}[l](-.71,-.71){$P_1$}
    \uput{\labelsep}[u](0,0){$P_2$}
    \uput{\labelsep}[u](1,0){$P_3$}
    \uput{\labelsep}[r](2,0){$P_4$}
    \uput{\labelsep}[r](2,1){$P_5$}}
\end{pspicture} 
\end{center}
and together with the weights $5,2,2,1,1$, it becomes the Enriques
diagram $\tT_{5,7}$ of the minimal log resolution of $x^5-y^7=0$.  In
general, the Enriques diagram $\tT_{p,q}$ that encodes the minimal log
resolution of the curve $x^p-y^q=0$ consists of the Enriques tree
$T_{p,q}$ together with the corresponding remainders of the Euclidean
algorithm as weights.
\end{exa}

\begin{cor}
  \label{c:star}
Let $D=\sum_\alpha e_\alpha E_\alpha$ be the divisor associated to the
Enriques diagram $\tT_{p,q}$.  If $P_\beta$ is a satellite point
proximate to $P_{\alpha'}$ and $P_{\alpha''}$, then
\[
  e_\beta = e_{\alpha'}+e_{\alpha''}+r_\beta.
\]  
\end{cor}

\proof
By Lemma \ref{l:star} and the interpretation of the weights of
$T_{p,q}$. 
\qed

\subsection{Multiplier ideals and jumping numbers}
\label{ss:multiplierIdeals}
We briefly recall the notions of multiplier ideals and jumping
numbers.  We refer the reader to \cite{La} for the results cited
below.  In the context of curves on surfaces, we define the relevant
divisors following \cite{SmTh}.

If $X$ is a smooth variety, $D\subset X$ an effective $\QQ$-divisor
and $\mu\colon Y\to X$ a {\it log resolution} for $D$, then
$\Jj(D)=\mu_\ast\Oo_Y(K_\mu-\floor{\mu^\ast D})$ is an ideal sheaf on
$X$%
\footnote{The integral part or round-down $\floor{D}$ of $D$ is the
integral divisor $\floor{D} = \sum_\alpha\floor{c_\alpha}D_\alpha$,
where for $c\in\QQ$, $\floor{x}$ denotes the greatest integer
$\leq c$.}.
The divisor $K_\mu$ is the relative canonical divisor of the map
$\mu$.  The ideal sheaf $\Jj(D)$ is independent of the choice of the
resolution and is called the multiplier ideal of $D$.  When $\mu$ is
the resolution of $D$ at a singular point $P$, the multiplier ideal
may be denoted by $\Jj(D)_P$.  The sheaf 
$\Oo_Y(K_\mu-\floor{\mu^\ast D})$ computing the multiplier ideal
satisfies, for any $i>0$, the local vanishing result
\[
  R^i\mu_\ast\Oo_Y(K_\mu-\floor{\mu^\ast D})=0.
\]

\begin{likethm}[Definition-Lemma]\rm (see \cite{EiLa})
Let $D\subset X$ be an effective divisor and $P\in D$ be a fixed
point.  Then there is an increasing discrete sequence of rational
numbers $\xi_i=\xi(D,P)$,
\[
  0 = \xi_0 < \xi_1 < \cdots
\]
such that $\Jj(\xi D)_P = \Jj(\xi_i D)_P$ for every
$\xi\in[\xi_i,\xi_{i+1})$, and
$\Jj(\xi_{i+1}D)_P\varsubsetneq\Jj(\xi_i D)_P$.  The rational numbers
$\xi_i$ are called {\it the jumping numbers} of $D$ at $P$.
\end{likethm}

The jumping numbers of $D$ at $P$ are periodic (see
\cite[Theorem~9.3.24]{La}) and that they are completely determined by
the ones that are less than $1$.  Therefore, in the sequel, we will
talk about the jumping numbers $<1$.

We have anticipated in the introduction that a set of candidates for
the jumping numbers is easy to provide in case $C$ is a curve singular
at $P$ on the smooth surface $X$.  Indeed, let $\mu:Y\to X$ be a log
resolution of $C$ with $\mu\inv(P)=\cup_{\alpha=1}^rE_r$.  Then
$K_\mu=\sum_{\alpha=1}^r W_\alpha=\sum_{\alpha=1}^rk_\alpha E_\alpha$,
with $k_\alpha>0$.  Writing 
$\mu^\ast C=\Ctilde+\sum_{\alpha=1}^r e_\alpha E_\alpha$, form the
proof of the above lemma it follows that the set of jumping numbers
must be contained in the set of the rational numbers
$(k_\alpha+n)/e_\alpha$, where $1\leq\alpha\leq r$ and $n$ is a
positive integer.

\begin{definition}[see \cite{SmTh}]
Let $\xi=(k_\alpha+n)/e_\alpha$ be a jumping number of $C$ at $P$.
The exceptional divisor $E_\alpha$ is said to contribute the jumping
number $\xi$ if
\[
  \Jj(\xi\cdot C) \varsubsetneq 
  \mu_\ast\Oo_Y(K_\mu-\floor{\xi\mu^\ast C}+E_\alpha).
\]
\end{definition}

If the above inclusion is satisfied for $\xi=(k_\alpha+n)/e_\alpha$,
then $\xi$ is a jumping number, since for any sufficiently small
$\epsilon>0$,
\[
  \mu_\ast\Oo_Y(K_\mu-\floor{\xi\mu^\ast C}+E_\alpha)
  \subset \mu_\ast\Oo_Y(K_\mu-\floor{(\xi-\epsilon)\mu^\ast C})
  = \Jj((\xi-\epsilon)\cdot C).
\]
In \cite{SmTh}, on the one hand, Theorem 3.1 shows that a divisor
$E_\alpha$ contributes to the sequence of jumping numbers if and only
if $E_\alpha$ has non-trivial intersection with at least three of the
other components of the total transform $\mu^\ast C$%
\footnote{The same characterization was established in 
\cite[Lemma~2.11]{FaJo} in an analytical context.}.  
On the other hand, it is also shown that not all jumping numbers are
contributed by exceptional divisors.  For example no $E_\alpha$
contributes the log-canonical threshold $1/2$ of the curve defined by
$(x^2-y^3)(x^3-y^2)=0$.  It is to be noticed that this curve is not
unibranch at the origin.  The essential step in the computation of the
jumping numbers of a unibranch curve is that each jumping number is
contributed by an exceptional divisor $E_\alpha$.  This is the content
of the forthcoming Proposition~\ref{p:unibranch}.

\begin{exa*}
In \cite{Ho} it is shown that if the
Puiseux exponent of $C$ at $P$ is $p/q$, with $\gcd(p,q)=1$, then the
jumping numbers less than $1$ of $C$ are $a/p+b/q<1$ with $a$ and $b$
positive integers.  There is only one divisor that contributes all
these jumping numbers, namely the last exceptional divisor.  
\end{exa*}

\begin{definition}
An exceptional divisor $E_\rho$ is said to be a {\it relevant
divisor}, or $\rho$ is said to be a {\it relevant position} of $C$ at
$P$, if $E_\rho\cdot E_\rho^0\geq3$, where 
$E_\rho^0=(\mu^\ast C)_{red}-E_\rho$.  The set of relevant positions
of $C$ at $P$ will be denoted by $\RRR_P$.
\end{definition}

\begin{rem}
  \label{r:branchBasis}
A relevant position $\rho$ is easy to identify on the Enriques tree.
Either it corresponds to a satellite point from which a 'slant'
edge starts, or it corresponds to a non-zero coefficient in the
expression of $\mu^\ast C-\Ctilde$ in the {\it branch basis}.  To
define this basis, 
let $\Pi$ be the proximity matrix.  The 
intersection matrix of the curves $W_\alpha$ is minus the
identity.  It follows that there exists effective%
\footnote{$||B_1\ldots B_s||= ^t\!\!\Pi\inv ||W_1\ldots W_s||$ and the
matrix $\Pi\inv$ has non negative entries since it decomposes the
$W_\alpha$ in terms of the $E_\alpha$.}
divisors $B_\alpha$ that form the dual basis to $(-E_\alpha)$ with
respect to the intersection form for the lattice
$\Lambda_C=\bigoplus_\alpha \ZZ E_\alpha$.  This basis is the
{\it branch basis}%
\footnote{The divisor $\mu^\ast C-\Ctilde$ may be expressed in three ways.  Its
expression in the branch basis reflects the branches of
$\Ctilde$---the analytically irreducible components of $\Ctilde$ above
$P$.  Its expression in the basis of total transforms $(W_\alpha)$
reflects the multiplicities of $C$ and the multiplicities of its
strict transforms along the resolution process, as it has been noticed
in Remark~\ref{r:theWeights}.  Finally, its expression in the basis
$(E_\alpha)$, the basis of strict transforms, gives the coefficients
necessary to compute the multiplier ideals associated to $C$.}.

Going back to the relevant position $\rho$ corresponding to a non-zero
coefficient in the branch basis, we notice that this position 
might be represented on the tree using arrowhead vertices, the number
of arrows being given by the coefficient $b_\rho$ in $\mu^\ast
C-\Ctilde=\sum_\alpha b_\alpha B_\alpha$.  We would obtain in
this way an augmented Enriques tree, equivalent to the Enriques
diagram $\tT$.
\end{rem}

\section{Jumping numbers of a unibranch curve}
\label{s:three}

The aim of this section is to prove Theorem~\ref{th:mainEnriques}.  To
formulate it we need to make some considerations about the {\it
unibranch trees}.  A unibranch tree is an Enriques tree having
out-valence $1$ for any of its vertices.  The trees $T_{p,q}$
introduced in Example~\ref{ex:tpq} are unibranch and represent the
simplest such trees in the sens that there is no slant edge starting
at a satellite point.  If a curve is unibranch at $P$, then the
Enriques diagram is given by a unibranch Enriques tree with the last
element in the branch basis $(B_\alpha)$ as the associated divisor.
The next definition allows us to see a unibranch tree as
being constructed from $T_{p,q}$ trees.

\begin{definition}
Let $T$ and $T'$ be unibranch Enriques trees with
$\VVV(T)=\{P_1,\ldots,P_r\}$ and $\VVV(T')=\{P'_1,\ldots,P'_{r'}\}$.
The {\it connected sum} of $T$ and $T'$ is the Enriques tree $T\#T'$
with the set of vertices
$\VVV(T\#T')=\VVV(T)\cup\VVV(T')/\{P_r=P'_1\}$, the set of edges
$\EEE(T\#T')=\EEE(T)\cup\EEE(T')$ and the map $\epsilon_{T\#T'}$
defined by $\epsilon_T$ and $\epsilon_{T'}$ through the natural
restrictions.
\end{definition}

\begin{exa}
  \label{ex:asUsual}
The minimal log resolution of $(x^3-y^2)^2-4x^5y-x^7=0$ needs five blowings
up with the following Enriques diagram.
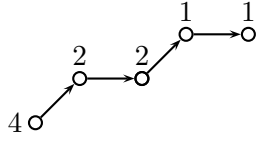
\begin{center}
\begin{pspicture}(0,0)(3,2.25)
  \psset{arrows=->,radius=.1,unit=5ex}
  \rput(1,1){
    \Cnode(-.71,-.71){1}
    \Cnode(0,0){2}
    \Cnode(1,0){3}
    \ncline{1}{2}
    \ncline{2}{3}
    \uput{\labelsep}[l](-.71,-.71){$4$}
    \uput{\labelsep}[u](0,0){$2$}
    \uput{\labelsep}[u](1,0){$2$}}
  \rput(2.71,1.71){
    \Cnode(-.71,-.71){1}
    \Cnode(0,0){2}
    \Cnode(1,0){3}
    \ncline{1}{2}
    \ncline{2}{3}
    \uput{\labelsep}[u](0,0){$1$}
    \uput{\labelsep}[u](1,0){$1$}}
\end{pspicture} 
\end{center}
The Enriques tree is the connected sum $T_{2,3}\#T_{2,3}$.
\end{exa}

\begin{thm}
  \label{th:mainEnriques}
Let $C$ be a curve unibranch at $P$ with the Enriques tree of the
minimal log resolution $S=T_{p_1,q_1}\#\cdots\#T_{p_g,q_g}$.  Set
$\qbar_1=q_1$ and
\[
  \qbar_j = \frac{m_{j-1}}{m_{j+1}}\,\qbar_{j-1}-p_j+q_j
\]
for any $2\leq j\leq g$, where $m_j = p_j\cdots p_g$ for any 
$1\leq j\leq g$, and $m_{g+1}=1$.  Then the jumping numbers less
that $1$ of $C$ at $P$ are given by
\[
  \bigcup_{j=1}^g \frac{1}{m_j \qbar_j}\, 
  R^{m_{j+1}}(p_j,\qbar_j),
\]
where
\[
  R^{m_{j+1}}(p_j,\qbar_j) 
  = \bigcup_{k=0}^{m_{j+1}-1}
  \big( kp_j\qbar_j+\{ap_j+b\qbar_j \mid a,b\in\NN^\ast,\,
  ap_j+b\qbar_j<p_j\qbar_j\} \big).
\]
\end{thm}

We begin by describing the sets $R^m(p,q)$ as they will appear in the
proof of Theorem~\ref{th:mainEnriques}.  For the purposes of this
section we denote by $\ceil{x}$ the round-up of $x$, \ie the least
integer $\geq x$, and by $\fpart{x}=x-\floor{x}$ the fractional part of
$x$.

Let $2\leq p<q$ be relatively prime integers and let $m$ be a positive
integer.  Setting $q'$ to be the positive integer that satisfies
$q'<p$ and $qq'=-1\mod p$, we define $R^m(p,q)$ as the set of integers
$k$, $1\leq k<mpq$, such that
\[
  \fpart{\frac{k}{pq}} + \fpart{\frac{q'k}{p}} > 1.
\]
If $m=1$ we shall denote the set $R^1(p,q)$ by $R(p,q)$.
Clearly 
\begin{equation} \label{eq:rInGeneral}
  R^m(p,q) = \bigcup_{j=0}^{m-1}(jpq+R(p,q))
\end{equation}
and $R(p,q)$ is computed in the following Proposition.

\begin{pro}
  \label{p:theSetR}
\[
  R(p,q) = \{ap+bq \mid a,b\in\NN^\ast,\, ap+bq<pq\}
\]
\end{pro}

\proof
If $k_0\in R(p,q)$ then 
\[
  \frac{k_0+p}{pq} + \fpart{\frac{q'(k_0+p)}{p}}
  = \frac{k_0}{pq}+\frac{1}{q} + \fpart{\frac{q'k_0}{p}}
  > 1,
\]
hence $k_0+p\in R(p,q)$.  It follows that to determine $R(p,q)$ it is
sufficient to determine the first element belonging to $R(p,q)$ in
each equivalence class mod $p$.  Clearly the multiples of $p$ do not
belong to $R(p,q)$.  So such an element is of the form $jq+Np$, with
$j\in\{1,\ldots,p-1\}$ and $N$ a positive integer.  Using the
hypothesis $qq'=-1\mod p$,
\[
  \frac{jq+Np}{pq}+\fpart{\frac{q'(jq+Np)}{p}} 
  = \frac{j}{p}+\frac{N}{q}+1-\frac{j}{p}
  = 1+\frac{N}{q}.
\]
So the minimal element in each equivalence class different from zero
is $jq+p$.  The result follows. 
\qed

As we have anticipated in \S\,\ref{s:one}, the computation of the
jumping numbers depends on the fact that each one of them is
contributed by an irreducible exceptional divisor.  More precisely we
have the following proposition whose proof will be given in
the last section.

\begin{pro}
  \label{p:unibranch}
Let $C$ be a unibranch curve at $P$ and let $\mu$ be a log resolution
such that $\mu^\ast C=\Ctilde+\sum_{\alpha\in\VVV}e_\alpha E_\alpha$.
If $\xi$ is a jumping number of $C$ at $P$, then there exists $\beta$
a relevant position such that $\floor{\xi e_\beta}=\xi e_\beta$ and
such that $E_\beta$ contributes $\xi$. 
\end{pro}

\likeproof[Proof of Theorem~\ref{th:mainEnriques}]
By Proposition~\ref{p:unibranch} each jumping number is contributed by
a relevant divisor and it is sufficient to compute the jumping numbers
contributed by each relevant divisor.  Remember that a relevant
divisor is an exceptional divisor that satisfies 
$E_\rho\cdot E_\rho^0\geq3$.

If $\rho$ is a relevant position, then tensoring the exact sequence 
\[
  0 \lra
  \Oo_Y \lra
  \Oo_Y(E_\rho) \lra
  \Oo_{E_\rho}(E_\rho|_{E_\rho}) \lra
  0
\]
with $\Oo_Y(K_\mu-\floor{\xi\mu^\ast C})$ and
pushing it down to $X$ give
\begin{multline*}
 0 \lra
  \mu_\ast\Oo_Y(K_\mu-\floor{\xi\mu^\ast C}) \lra
  \mu_\ast\Oo_Y(K_\mu-\floor{\xi\mu^\ast C}+E_\rho) \\
  \lra \mu_\ast\Oo_{E_\rho}(K_{E_\rho}-\floor{\xi\mu^\ast C}|_{E_\rho})
  \lra 0
\end{multline*}
thanks to the local vanishing.  Since $E_\rho$ is a projective line,
$\xi$ is a jumping number contributed by $E_\rho$ if and only if 
$\xi e_\rho$ is an integer and
\begin{equation}
  \label{eq:mainInequality}
  -\floor{\xi\mu^\ast C} \cdot E_\rho \geq 2.
\end{equation}
Assume that $\xi e_\rho$ is an integer.  Let $r_j$ be the number of
vertices of the Enriques tree $T_{p_j,q_j}$ and set
$s=r_1+\cdots+r_g-(g-1)$, the number of vertices of $S$.  There are
two cases to be considered: either $\rho=s$, \ie $\rho$ is the highest
point of the Enriques tree, or $\rho$ is a relevant position different
from $s$.  Whatever the case, the study of the numbers
$\floor{\xi\mu^\ast C} \cdot E_\rho$ depends on the control of the
coefficients $e_\alpha$ in 
$\mu^\ast C=\Ctilde+\sum_{\alpha=1}^s e_\alpha E_\alpha$.  The
following technical lemmas \ref{l:eCoeff}, \ref{l:eCoeffBefore_r},
\ref{l:eRelevantCoeff} and \ref{l:eCoeffAfter_r}
give explicit formulae for these
coefficients and also, some other useful numerical relations.

Some notation is in order.  Let $T$ be an Enriques tree.  Denote by
$E_\alpha^T$ the elements of the basis of strict transforms, by
$W_\alpha^T$ the elements of the basis of total transforms and by
$B_\alpha^T$ the elements of the branch basis that has been introduced
in Remark~\ref{r:branchBasis}.  If 
$\Lambda=\bigoplus_\alpha E_\alpha^T$, then $(e_\alpha^T)$ will denote
the basis for $\Lambda^\ast$, dual to $(E_\alpha^T)$.  Similarly,
$(w_\alpha^T)$ will denote the dual basis to the basis of total
transforms and $(b_\alpha^T)$ the dual basis to the branch basis.

Let $p<q$ be two relatively prime positive integers.  Consider the
Euclidean algorithm $r_0=a_1r_1+r_2$,
\ldots, $r_{m-2}=a_{m-1}r_{m-1}+r_{m}$ and $r_{m-1}=a_mr_m$, with
$r_0=q$ and $r_1=p$.  Define as in \cite[Lemma~A.8]{Na} two finite
sequences 
 $(f_j)_{-1\leq j\leq m}$ and 
$(\delta_j)_{1\leq j\leq m+1}$ by
\begin{equation}
  \label{eq:sequences}
\begin{aligned}
  f_j &= f_{j-2}+a_j\delta_j, \quad\text{for any } 1\leq j\leq m, \\
  \delta_j &= \delta_{j-2}+a_{j-1}\,f_{j-2}, 
  \quad\text{for any } 2\leq j\leq m+1, \quad 
\end{aligned}
\end{equation}
such that $f_{-1}=f_0=0$ and $\delta_0=\delta_1=1$.  It is easy to show
that the remainder $r_j$ in the Euclidean algorithm is given by
$-f_{j-1}q+\delta_jp$ if $j$ is odd and by $\delta_jq-f_{j-1}p$ if $j$
is even.  Furthermore, if $m$ is odd, then $f_m=q$ and
$\delta_{m+1}=p$, and if $m$ is even, then $f_m=p$ and
$\delta_{m+1}=q$.  We have the following lemma that computes various
coefficients for the Enriques tree $T_{p,q}$.

\begin{lem}
  \label{l:eCoeff}
Let $T=T_{p,q}$.  Then for any $0\leq j\leq m-1$ and any 
$1\leq k\leq a_{j+1}$, 
\[
  e_{a_1+\cdots+a_j+k}^T(B_r^T) =
  \begin{cases}
    (f_{j-1}+k\delta_{j+1})\,p & \textql{if $j$ is even}\\
    (f_{j-1}+k\delta_{j+1})\,q & \textql{if $j$ is odd}
  \end{cases}
\]
and
\[
  e_{a_1+\cdots+a_j+k}^T(W_1^T) =
  \begin{cases}
    \delta_j+kf_j & \textql{if $j$ is even}\\
    f_{j-1}+k\delta_{j+1} & \textql{if $j$ is odd}.
  \end{cases}
\]
\end{lem}

\proof
We proceed by
induction using Corollary~\ref{c:star}, the relations
(\ref{eq:sequences}) and the relations quoted after them.  For the
computation of $e_{a_1+\cdots+a_j+k}^T(B_r^T)$, we suppose
that $j$ is even, the case $j$ odd being similar.  If $k=1$, then
\[
\begin{split}
  e_{a_1+\cdots+a_j+1}^T(B_r^T)
  &= e_{a_1+\cdots+a_{j-1}}^T(B_r^T) + e_{a_1+\cdots+a_j}^T(B_r^T)
  + r_{j+1} \\
  &= f_{j-1}p+f_jq+(-f_jq+\delta_{j+1}p) \\
  &= (f_{j-1}+\delta_{j+1})\,p.
\end{split}
\]
Now, if $1<k\leq a_{j+1}$, then
\[
\begin{split}
  e_{a_1+\cdots+a_j+k}^T(B_r^T)
  &= e_{a_1+\cdots+a_{j}}^T(B_r^T) + e_{a_1+\cdots+a_j+k-1}^T(B_r^T)
  + r_{j+1} \\
  &= f_{j}q + (f_{j-1}+(k-1)\delta_{j+1})\,p + (-f_jq+\delta_{j+1}p) \\
  &= (f_{j-1}+k\delta_{j+1})\,p.
\end{split}
\]
As for the second equality, if we suppose again that $j$ is even,
we get
\[
  e_{a_1+\cdots+a_j+1}^T(W_1^T)
  = e_{a_1+\cdots+a_{j-1}}^T(W_1^T) + e_{a_1+\cdots+a_j}^T(W_1^T)
  = \delta_j+f_j
\]
for $k=1$, and
\[
  e_{a_1+\cdots+a_j+k}^T(W_1^T)
  = e_{a_1+\cdots+a_{j}}^T(W_1^T) + e_{a_1+\cdots+a_j+k-1}^T(W_1^T)
  = f_j+(\delta_j+(k-1)f_j)
\]
for $1<k\leq a_{j+1}$.
\qed

The remaining lemmas will allow us to compute the $e_\alpha$
coefficients for $\mu^\ast C-\Ctilde$ for a unibranch curve.

\begin{lem}
  \label{l:eCoeffBefore_r}
Let $S=T\# T'$ be a unibranch Enriques tree.  If $r$ is the number of
vertices of $T$ and $s$ is the number of vertices of $S$, then
\[
  e_{\alpha}^S(B_{s}^S) =
  w_r^S(B_s^S)\,e_\alpha^T(B_r^T)
\]
for any $1\leq \alpha\leq r$.  
\end{lem}

\proof
Write $B_{s}^S=\sum_\beta w_\beta^S(B_{s}^S)\,W_\beta^S$ using the
basis of total transforms.  Since the proximity matrix in upper
unitriangular, if $1\leq \alpha\leq r$, then only the divisors
$W_\beta$ with $\beta\leq r$ count when computing
$e_{\alpha}^S(B_{s}^S)$.  For $\beta\leq r$, we have
\[  
  w_\beta^S(B_{s}^S) = w_\beta^S( w_r^S(B_{s}^{S})\,B_r^S)
  = w_r^S(B_s^S)\,w_\beta^S(B_r^S) = w_r^S(B_s^S)\,w_\beta^T(B_r^T).
\]
So
\[
\begin{split}
  e_\alpha^S(B_s^S) 
  &= w_r^S(B_s^S)\sum_{\beta\leq
    r}w_\beta^T(B_r^T)e_\alpha^S(W_\beta^S)\\
  &= w_r^S(B_s^S)\,e_\alpha^T\big(
  \sum_{\beta\leq r}w_\beta^T(B_r^T)W_\beta^T\big),
\end{split}
\]
since $e_\alpha^S(W_\beta^S)=e_\alpha^T(W_\beta^T)$, and the result
follows. 
\qed

\begin{lem}
  \label{l:eRelevantCoeff}
Let 
\[
  S = T_{p_1,q_1} \# T_{p_2,q_2} \#\cdots\# T_{p_g,q_g}
\]
with $r_j$ the number of vertices of the tree $T_{p_j,q_j}$, and let
$s=r_1+\cdots+r_g-(g-1)$ be the number of vertices of $S$.  
If $\qbar_1=q_1$ and
\[
  \qbar_j = p_{j-1}p_j\qbar_{j-1}-p_j+q_j
\]
for any $2\leq j\leq g$, then
\[
\begin{aligned}
  w_{r_1+\cdots+r_j-(j-1)}^S(B_s^S) &= p_j\cdots p_g\\
  e_{r_1+\cdots+r_j-(j-1)}^S(B_s^S) &= p_j\cdots p_g\qbar_j 
\end{aligned}
\]
for any $1\leq j\leq g$.
\end{lem}

\proof
Set $T_j=T_{p_j,q_j}$ for any $j$.  The first identity is clear by
induction since using the proximity relations,
\[
  w_{r_1+\cdots+r_j-(j-1)}^S(B_s^S) = 
  w_1^{T_{j+1}\#\cdots\# T_g}
  (B_{r_{j+1}+\cdots+r_g-(g-j-1)}^{T_{j+1}\#\cdots\# T_g}).
\]
As for the second, by Lemma~\ref{l:eCoeffBefore_r},
\[
\begin{split}
  e_{r_1+\cdots+r_j-(j-1)}^S(B_s^S) 
  &= w_{r_1+\cdots+r_j-(j-1)}^S(B_s^S)\, 
  e_{r_1+\cdots+r_j-(j-1)}^{T_1\#\cdots\# T_j}
    (B_{r_1+\cdots+r_j-(j-1)}^{T_1\#\cdots\# T_j}) \\
  &= p_{j+1}\cdots p_gp_j\qbar_j. 
\end{split}
\]
\qed

\begin{lem}
  \label{l:eCoeffAfter_r}
Let $S=T'\# T$ be a unibranch Enriques tree with $T=T_{p,q}$.  Let
$r'$ be the number of vertices of $T'$ and $r$ the number of vertices
of $T$.  Then
\[
  e_{r'-1+\alpha}^S(B_{r'-1+r}^S)
  = \big( e_{r'}^{T'}(B_{r'}^{T'})-1\big)\,
  e_\alpha^T(W_1^T)\, p + e_\alpha^T(B_r^T)
\]
for any $1\leq \alpha\leq r$.
\end{lem}

\proof 
Set $s=r'+r-1$, the number of vertices of $S$.  
Using that
\[
  w_{r'-1+\beta}^{S}(B_{s}^{S})=w_{\beta}^T(B_r^T) 
\]
for any $2\leq\beta\leq r$ and discarding the exponent $S$, we have
\[
\begin{split}
  e_{r'-1+\alpha}(B_{s}) &=
  e_{r'-1+\alpha}\bigg(\sum_{\gamma=1}^{r'}w_\gamma(B_{s})\,W_\gamma
  + \sum_{\beta=2}^rw_{r'-1+\beta}(B_{s})\,W_{r'-1+\beta}
  \bigg)\\
  &= e_{r'-1+\alpha}\bigg(
  w_{r'}(B_{s})\,B_{r'} - w_{r'}(B_{s})\,W_{r'}
  + \sum_{\beta=1}^rw_{r'-1+\beta}(B_{s})\,W_{r'-1+\beta}
  \bigg)\\
  &= e_{r'-1+\alpha}
  \big( w_1^T(B_r^T)\,B_{r'} - w_1^T(B_r^T)\,W_{r'} \big)
  + e_\alpha^T\bigg(
  \sum_{\beta=1}^rw_\beta^T(B_r^T)\,W_\beta^T
  \bigg)\\
  &=  w_1^T(B_r^T)\,
  e_{r'-1+\alpha}\big(e_{r'}^{T'}(B_{r'}^{T'})W_{r'}-W_{r'}\big)
  + e_\alpha^T(B_r^T),
\end{split}
\]
and hence the result, since 
$e_{r'-1+\alpha}(W_{r'})= e_\alpha^T(W_1^T)$ and $w_1^T(B_r^T)=p$.
\qed

\likeproof[End of the proof of Theorem~\ref{th:mainEnriques}]
We need to study $\floor{\xi\mu^\ast C} \cdot E_\rho$ when $\rho$ is a
relevant position and $\xi e_\rho$ is an integer.  We have already
noticed that there are two cases to be considered: either $\rho=s$ or
$\rho\neq s$, where $s$ is the number of vertices of
$S=T_{p_1,q_1}\#\cdots\#T_{p_g,q_g}$.

In the former case, set $T=T_{p_g,q_g}$.  Then
\[
  \floor{\xi\mu^\ast C} \cdot E_{s}^S
  = \Big(
  \floor{\xi e_{s-a_m}}E_{r_g-a_m}^{T}+
  \floor{\xi e_{s-1}}E_{r_g-1}^{T}+
  \xi e_{s}E_{r_g}^{T}
  \Big) \cdot E_{r_g}^{T},
\]
where $a_m$ is the last quotient in the Euclidean algorithm for $q_g$
and $p_g$. By Corollary~\ref{c:star}, $e_{s}=e_{s-a_m}+e_{s-1}+1$.
By Lemma~\ref{l:eRelevantCoeff}, $e_{s}=p_g\qbar_g$ with
$p_g<\qbar_g$, relatively prime integers.
Then, by Lemma~\ref{l:eCoeffAfter_r}
\[
  e_{s-\alpha}
  = (p_{g-1}\qbar_{g-1}-1)\, p_g\, e_{r_g-\alpha}^T(W_1^T)
  + e_{r_g-\alpha}^T(B_{r_g}^T)
\]
for any $1\leq\alpha<r_g$.  One of the two positions $r_g-m$ and
$r_g-1$ must belong to a proper horizontal $L$-shape branch.  We
suppose that $r_g-1$ satisfies this, the argument being similar in the
other case.  By Lemma~\ref{l:eCoeff}, 
$e_{r_g-1}^T(B_{r_g}^T) = e_{r_g-1}^T(W_1^T)\, q_g$, \ie
\[
  e_{s-1} = (p_{g-1}\qbar_{g-1}-1)\, p_g\, e_{r_g-1}^T(W_1^T)
  + e_{r_g-1}^T(B_{r_g}^T)
  = e_{r_g-1}^T(W_1^T)\, \qbar_g,
\]
and $p_g$ divides $e_{s-m}$.  Set $M=e_{r_g-1}^T(W_1^T)$.  From
$p_g\qbar_g = e_{s-m}+M\qbar_g+1$
it follows that 
\[
  M\qbar_g = -1 \mod p_g.
\]
Set $x=\xi e_{s}=\xi p_g\qbar_g$.  It is an integer
satisfying $1\leq x<p_g\qbar_g$.  Then putting everything together,
\[
\begin{split}
  -\floor{\xi\mu^\ast C} \cdot E_{s}
  &=
  -\floor{\frac{x}{p_g\qbar_g}(p_g\qbar_g-M\qbar_g-1)}
  -\floor{\frac{x}{p_g\qbar_g}\,M\qbar_g}+x\\
  &=\ceil{\frac{Mx}{p_g}+\frac{x}{p_g\qbar_g}}-\floor{\frac{Mx}{p_g}}.
\end{split}
\]
Hence the inequality (\ref{eq:mainInequality}) is satisfied, \ie
$-\floor{\xi\mu^\ast C} \cdot E_\rho \geq 2$, if and only if
\[
  \fpart{\frac{Mx}{p_g}}+\frac{x}{p_g\qbar_g}>1
\]
By Proposition~\ref{p:theSetR} this is equivalent to 
$x\in R(p_g,\qbar_g)$.

In the second case, if $\rho=r_1+\cdots+r_j-(j-1)$ is a relevant
position different from $s$, the highest one, set $r=r_{j}$,
$T=T_{p_{j},q_{j}}$ and $S=T'\# T\# T''$.  When computing
$\floor{\xi\mu^\ast C}\cdot E_\rho$ we distinguish two situations
according to whether $a''_1=1$ or not, where $a''_1,a''_2\ldots$ are
the quotients in the Euclidean algorithm for $q_{j+1}$ and $p_{j+1}$.
In case $a''_1\neq1$,
\[
\begin{split}
  \floor{\xi\mu^\ast C} \cdot E_\rho
  &= \Big(
  \floor{\xi e_{\rho-a_m}}E_{\rho-a_m}+
  \floor{\xi e_{\rho-1}}E_{\rho-1}+
  \xi e_\rho E_\rho+
  \floor{\xi e_{\rho+1}}E_{\rho+1}
  \Big) \cdot E_\rho \\
  &= \floor{\xi e_{\rho-a_m}}+\floor{\xi e_{\rho-1}}-2\xi e_\rho
  +\floor{\xi e_{\rho+1}},
\end{split}
\]
whereas in case $a''_1=1$,
\[
  \floor{\xi\mu^\ast C} \cdot E_\rho=
  \floor{\xi e_{\rho-a_m}}+\floor{\xi e_{\rho-1}}-(2+a''_2)\xi e_\rho
  +\floor{\xi e_{\rho+a''_2+1}}.
\]
By Lemma~\ref{l:star}, $e_\rho=e_{\rho-a_m}+e_{\rho-1}+m_{j+1}$, and
by Lemma~\ref{l:eRelevantCoeff}, $e_\rho=p_j\qbar_jm_{j+1}$.  As
before, either $m$ is odd and $e_{\rho-a_m}=M\qbar_jm_{j+1}$, or $m$
is even and $e_{\rho-1}=M\qbar_jm_{j+1}$ with 
$M\qbar_j = -1 \mod p_j$.  Set $x=\xi e_\rho$.  The integer $x$
satisfies $1\leq x<e_\rho=p_j\qbar_jm_{j+1}$.  We claim that
independently of $a''_1$,
\begin{equation} \label{eq:claimStep1}
  -\floor{\xi\mu^\ast C} \cdot E_\rho = 
  \ceil{\frac{Mx}{p_j}+\frac{x}{p_j\qbar_j}}
  -\floor{\frac{Mx}{p_j}}-\floor{\frac{x}{p_j\qbar_j}}.
\end{equation}
The justification is more complicated if $a''_1=1$.  The Enriques tree
is shown below in case $a''_1=1$ and $m$ even---the last proper
$L$-shape branch of $T_{p_j,q_j}$ is horizontal.
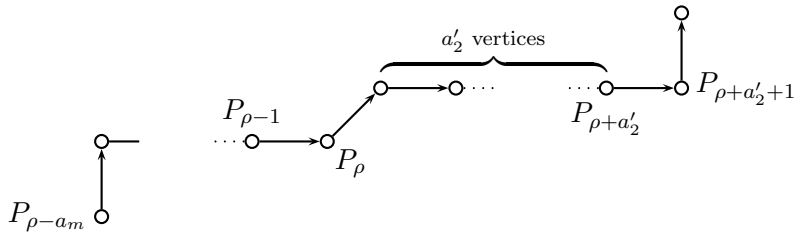
\begin{center}
  \begin{pspicture}(0,0)(12,3.8)
    \rput(5.5,2.3){
      \uput{\labelsep}[u](1.5,0){$\substack{a'_2\text{ vertices}
          \\ \overbrace{\hspace{30mm}}}$}
      \psset{arrows=->,radius=.1}  
      \Cnode(-3.71,-1.71){13}
      \Cnode(-3.71,-.71){12}
      \Cnode(-1.71,-.71){11}
      \Cnode(-.71,-.71){1}
      \ncline{13}{12}
      \ncline{11}{1}
      \Cnode(0,0){2}
      \Cnode(1,0){3}
      \Cnode(3,0){4}
      \Cnode(4,0){5}
      \Cnode(4,1){6}
      \ncline{1}{2}
      \ncline{2}{3}
      \ncline{4}{5}
      \ncline{5}{6}
      \uput{\labelsep}[l](-3.71,-1.71){$P_{\rho-a_m}$}
      \uput{\labelsep}[u](-1.71,-.71){$P_{\rho-1}$}
      \uput{\labelsep}[l](0,-1){$P_\rho$}
      \uput{\labelsep}[d](3,0){$P_{\rho+a'_2}$}
      \uput{\labelsep}[r](4,0){$P_{\rho+a'_2+1}$}
      \psset{arrows=-}
      \psline(-3.61,-.71)(-3.21,-.71)
      \psline[linestyle=dotted](-2.21,-.71)(-1.81,-.71)
      \psline[linestyle=dotted](1.1,0)(1.5,0)
      \psline[linestyle=dotted](2.5,0)(2.9,0)}
  \end{pspicture} 
\end{center}
By Lemma~\ref{l:eCoeff},
\[
  e_{\rho+a''_2+1} = (a''_2+1)e_\rho + a''_2r''_2m_{j+2} + r''_3m_{j+2} 
  = (a'_2+1)p_j\qbar_jm_{j+1}+m_{j+1},
\]
hence $-\floor{\xi\mu^\ast C}\cdot E_\rho$ is given by 
\begin{multline*}
  -\floor{\frac{x}{p_j\qbar_jm_{j+1}}(p_j\qbar_jm_{j+1}
    -M\qbar_jm_{j+1}-m_{j+1})}
  -\floor{\frac{x}{p_j\qbar_jm_{j+1}}\,M\qbar_jm_{j+1}} \\
  +(2+a'_2)x
  -\floor{\frac{x}{p_j\qbar_jm_{j+1}}((a'_2+1)p_j\qbar_jm_{j+1}
    +m_{j+1})}  
\end{multline*}
and the equality follows establishing the equality
(\ref{eq:claimStep1}).  Finally, 
$-\floor{\xi\mu^\ast C} \cdot E_\rho \geq 2$ is equivalent to
\[
  \fpart{\frac{Mx}{p_j}}+\fpart{\frac{x}{p_j\qbar_j}}>1
\]
with $x=\xi e_\rho<p_j\qbar_jm_{j+1}$, \ie to 
$x\in R^{m_{j+1}}(p_j,\qbar_j)$, again by Proposition~\ref{p:theSetR}.
\qed

\begin{exa*}
A jumping number might be contributed by more than one relevant
divisor.  For example if the Enriques diagram associated to the
minimal log resolution of $C$ is given by 
$T = T_{2,3}\# T_{5,11}$---a tree with nine vertices---and
$\mu^\ast C = \Ctilde+B_9^T$, then $\xi=11/30$ is a jumping number
contributed by either $E_3$ or $E_9$.
\end{exa*}

The first jumping numbers of a unibranch curve can be obtained by
inspecting the jumping numbers of the term ideal associated to $C$ in
a suitable coordinate system.

\begin{cor}
Let $C$ be a curve unibranch at $P$ with the Enriques tree of the
mi\-ni\-mal log resolution $T_{p_1,q_1}\#\cdots\#T_{p_g,q_g}$.  Fixing
$\Pp$ an allowable system of local parameters, let $\aaa_{C,\Pp}$
be the term ideal of an equation of $C$ in $\Pp$.  Then the first 
$\card(R(p_1,q_1))$ jumping numbers of $C$ at $P$ coincide with the first
jumping numbers of $\aaa_{C,\Pp}$.
\end{cor}

\proof
Set $\pi=p_2\cdots p_g$.  If $\tT$ is the Enriques diagram defined by
the tree $T_{p_1,q_1}$ and the wights corresponding to the last
$B_\alpha$ considered with multiplicity $\pi$, then 
\[
  \overline{\aaa}_{C,\Pp} = \mu^\ast\Oo_Y(-D_\tT).
\] 
Notice that $\overline{\aaa}_{C,\Pp}$ is the smallest integrally
closed monomial ideal containing an equation of $C$.  The jumping
numbers of $\aaa_{C,\Pp}$, or equivalently the jumping numbers of 
$\overline{\aaa}_{C,\Pp}$, less than $1$ are given by 
\[
  \frac{ap_1+bq_1}{p_1q_1\pi} < 1,
\]
with $a,b$ positive integers.  Then the first $\card(R(p_1,q_1))$ of
them, \ie those for which $ap_1+bq_1<p_1q_1$ are among the jumping
numbers of $C$, in the subset $R^\pi(p_1,q_1)$, by
Theorem~\ref{th:main}.  It is sufficient to show that these jumping
numbers are also the first ones of $C$.  But they all satisfy
\[
  \frac{ap_1+bq_1}{p_1q_1\pi} < \frac{1}{\pi}
\] 
and $1/\pi$ is bigger that any element of any set
$1/m_j\qbar_jR^{m_{j+1}}(p_j,\qbar_j)$, with $j\geq 2$.
\qed

\section{Reformulation of Theorem~\ref{th:mainEnriques} in terms of \\the
  semigroup of the singularity}
\label{s:two}

Let $S(C)$ be the semigroup of the curve $C$ unibranch at $P$.  It is
defined by
\[
  S(C) = \{\ord_P s \mid s\in\Oo_{C,P}\}
\]
where the order of the local section $s$ is computed using a
normalization $\Ctilde\to C$.  If the Puiseux characteristic of $C$ at
$P$ is $(m;\beta_1,\ldots,\beta_g)$, the first part of Theorem~4.3.5
in \cite{Wa} states that the integers $\betabar_0$,
$\betabar_1,\ldots,\betabar_g$ generate $S(C)$, where $\betabar_j$ are
defined by $\betabar_0=m=m_1$, $\betabar_1=\beta_1$ and then
inductively by
\begin{equation}
  \label{eq:betabar}
  \betabar_{j+1} = \frac{m_j}{m_{j+1}}\,\betabar_j+\beta_{j+1}
  - \beta_j,
\end{equation}
for any $1\leq j<g$, and where $m_{j+1}=\gcd(m_j,\beta_j)$ for any
$1\leq j\leq g$.  The second part shows how to recover the generators
$\betabar_j$ once the semigroup $S(C)$ given: $\betabar_0=m_1$ is
the least non-zero element of $S(C)$, and inductively $\betabar_{j}$
is the least element non divisible by $m_j$, and
$m_{j+1}=\gcd(m_j,\betabar_j)$.  We can state the main result of the
paper.

\begin{thm}
  \label{th:main}
Let $C$ be a curve unibranch at $P$ and let
$(\betabar_0,\betabar_1,\ldots,\betabar_g)$ be the canonical minimal
set of generators of the semigroup $S(C)$.  Then the jumping numbers
of $C$ at $P$ less that $1$ are given by
\[
  \bigcup_{j=1}^g \frac{1}{[m_j,\betabar_j]}\, 
  R^{m_{j+1}}\bigg(\frac{m_j}{m_{j+1}},\frac{\betabar_j}{m_{j+1}}\bigg),
\]
where $m_1=\betabar_0$ and $m_{j+1}=\gcd(m_j,\betabar_j)$ for any
$j\geq1$, and where $[m_j,\betabar_j]$ denotes the least common
multiple of the two integers.
\end{thm}

For the proof we will use Theorem~\ref{th:mainEnriques} and need
Enriques' equivalence between the Puiseux characteristic of $C$ at
$P$ and the Enriques diagram associated to the minimal log resolution
of $C$ at $P$, equivalence that we present next.

Let $(x,y)$, be a system of local parameters.  If $x=0$ is not tangent
to $C$ at $P$, there exists a good parametrization
(see~\cite[Chapter~2]{Wa}) of $C$, $x=t^m$ and 
$y=\sum_{k=m}^\infty c_kt^k$.  The {\it Puiseux characteristic} of $C$
is the sequence of integers $(m;\beta_1,\ldots,\beta_g)$ defined as
follows: $\beta_1$ is the exponent of the first term in the power
series which is not a power of $t^m$.  Set $m_1=m$ and
$m_2=\gcd(m_1,\beta_1)$.  Inductively, $\beta_j$ is the exponent of
the first term which is not a power of $m_j$ and
$m_{j+1}=\gcd(m_j,\beta_j)$.  The construction stops when $m_{g+1}=1$
is reached.  The integers $\beta_j$ are also called the Puiseux
characteristic exponents and the {\it Puiseux exponents} are just the
rationals $\beta_1/m, \beta_2/m,\ldots,\beta_g/m$.  Note that
\[
  \frac{\beta_1}{m} < \frac{\beta_2}{m} < \cdots < 
  \frac{\beta_g}{m}.
\]

\begin{pro}[see \cite{Wa}, Theorem 3.5.5]
  \label{p:puiseuxDropping}
Let $C\subset X$ be a unibranch curve at $P$ on the smooth surface
$X$.  If the Puiseux characteristic of $C$ at $P$ is
$(m;\beta_1,\ldots,\beta_g)$, then the Puiseux characteristic of the
curve obtained by blowing up $X$ at $P$ is given by
\[
\begin{split}
  (m;\beta_1-m,\ldots,\beta_g-m) &\textq{if} \beta_1>2m\\
  (\beta_1-m;m,\beta_2-\beta_1+m,\ldots,\beta_g-\beta_1+m)
  &\textq{if} \beta_1<2m, \, (\beta_1-m)\!\!\not|m\\
  (\beta_1-m;\beta_2-\beta_1+m,\ldots,\beta_g-\beta_1+m)
  &\textq{if} (\beta_1-m)|m.
\end{split}
\]
\end{pro}

\begin{cor}
  \label{c:puiseuxEnriques}
Let $C$ be a curve unibranch at $P$ whose Puiseux characteristic is\linebreak
$(m;\beta_1,\ldots,\beta_g)$.  Consider the sequence of the strict
transforms of $C$ constructed by the successive blowings up that
desingularize $C$.  The first Puiseux characteristic in this sequence
with exactly $g-j+1$ characteristic exponents is
$(m_j;\beta_j^{(j-1)},\ldots,\beta_g^{(j-1)})$, where 
$m_j= \gcd(m_{j-1},\beta_{j-1})$ as before, and for any $k\geq j$,
\[
  \beta_k^{(j-1)} = \beta_k-\beta_{j-1}+a_{j-1}m_j
\]
with $a_{j-1}$ the last quotient in the Euclidean
algorithm for $m_{j-1}$ and $\beta_{j-1}$. 
\end{cor}

\proof
Considering the Euclidean algorithm for $\beta_1$ and $m_1=m$ and
using the previous proposition, it is easy to see that the first
element in the second Puiseux characteristic is
$m_2=\gcd(m_1,\beta_1)$.  Computing $\beta_j-\beta'_j$ we get that
that $\beta'_j=\beta_j-\beta_1+a_1m_2$, for any $j\geq2$.  Notice that
\[
  \gcd(m_2,\beta_2) = \gcd(m_2,\beta'_2)
\]
and that furthermore, the last quotients in the Euclidean algorithms
for $\beta_2$ and $m_2$, and $\beta'_2$ and $m_2$ coincide.
The result follows by induction on $j$.
\qed

Enriques established the equivalence between the Puiseux
characteristic and the Enriques diagram associated to a curve
unibranch and singular at $P$.  We refer the reader to \cite{EnCh},
\cite{Ca} and especially \cite{SeKn} Theorem XI~\S~6.1.3 for further
details.  What we want to state now is just a more concise way to
present this equivalence.

\begin{thm}[Enriques]
  \label{th:enriques}
Let $(m;\beta_1,\ldots,\beta_g)$ be the Puiseux characteristic of $C$
at $P$.  Set 
\[
  p_j = \frac{m_j}{m_{j+1}} 
  \textq{and}
  q_j = \frac{\beta_j-\beta_{j-1}+m_j}{m_{j+1}} 
\]
for any $1\leq j\leq g$, with $\beta_0=m_1$.  Then, the corresponding
Enriques diagram $\tT$ is given by the Enriques tree
\[
  T = T_{p_1,q_1} \# T_{p_2,q_2} \#\cdots\# T_{p_g,q_g}.
\]
\end{thm}

\likeproof[Sketch of proof]\!\!\!
We will use the notation from Corollary~\ref{c:puiseuxEnriques}.  The
first part of the Enriques tree associated to $C$ coincides with the
Enriques tree associated to a curve having the Puiseux characteristic
$(m_1;\beta_1)$.  Such a curve is desingularized by the blowings up
corresponding to the whole Enriques tree
$T_{p_1,\beta_1/m_2}=T_{p_1,q_1}$ except for the last stretch.  It is
noteworthy that the length of this stretch equals the last quotient in
the Euclidean algorithm for $\beta_1/m_2$ and $p_1$, \ie $a_1$, and
that the blowings up of this last stretch are needed to obtain a
log-resolution for such a curve.  Now, if $C_2$ denotes the strict
transform of $C$ before the blowings up of the last stretch of
$T_{p_1,\beta_1/m_2}$, then $C_2$ is the first strict transform among
the strict transforms of $C$ given by the log resolution, having
exactly $g-1$ Puiseux exponents.  The Puiseux characteristic is
$(m_2;\beta'_2,\ldots,\beta'_g)$.  What we have previously said for
$C$ is true for $C_2$.  So the Enriques tree associated to $C_2$
starts with $T_{p_2,n_2}$, where
\[
  n_2 = \beta'_2/m_3.
\]
It follows that the part of the Enriques tree associated to $C$ and
corresponding to the first two Puiseux exponents is 
\[
  T_{p_1,q_1} \# T_{p_2,n_2-(a_1-1)p_2}.
\]
But
\[
  n_2-(a_1-1)p_2 = \frac{\beta'_2}{m_3}-(a_1-1)\frac{m_2}{m_3}
  = \frac{\beta_2-\beta_1+a_1m_2}{m_3}-(a_1-1)\frac{m_2}{m_3}
  = q_2.
\]
The proof is finished by induction on the number of Puiseux exponents.
\qed

\begin{exa*}
The curve given in Example~\ref{ex:asUsual} has the good
parametrization $t\mapsto (t^4,t^6+t^7)$.  It follows that the Puiseux
characteristic is $(4;6,7)$, with $3/2$ the first Puiseux exponent,
and after one blow-up, the strict transform has the unique Puiseux
exponent $5/2$.  We recover the fact that the Enriques tree of the
minimal log resolution is $T_{2,3}\# T_{2,3}$.
\end{exa*}

\begin{cor}
  \label{c:wCoeff}
Let $C$ be a curve unibranch at $P$.  In the notation of
Theorem~\ref{th:enriques}, if $w_\alpha$ are the weights of the
Enriques diagram associated to the minimal log resolution of $C$ at
$P$, and if $r_j$ is the number of vertices of the tree $T_{p_j,q_j}$
for any $1\leq j\leq g$, then $w_1=m_1$ and 
\[
  w_{r_1+\cdots+r_j-(j-1)} = m_{j+1}
  = p_{j+1}p_{j+2}\cdots p_{g}
\]
for any $1\leq j\leq g$.  
\end{cor}

\proof
Follows from Theorem~\ref{th:enriques},
Proposition~\ref{p:puiseuxDropping}, Corollary~\ref{c:puiseuxEnriques}
and Remark~\ref{r:theWeights}.
\qed

\likeproof[Proof of Theorem~\ref{th:main}]
By Theorem~\ref{th:mainEnriques}, the jumping numbers are given by
\[
  \bigcup_{j=1}^g \frac{1}{m_j \qbar_j}\, 
  R^{m_{j+1}}(p_j,\qbar_j).
\]
where the Enriques tree of the minimal log resolution is
$T_{p_1,q_1}\#\cdots\# T_{p_g,q_g}$.  Furthermore
$m_j = p_j\cdots p_g$ and
\[
  \qbar_j = \frac{m_{j-1}}{m_{j+1}}\,\qbar_{j-1}-p_j+q_j
\]
for any $1\leq j\leq g$, with $m_{g+1}=1$ and $\qbar_1=q_1$.
Clearly 
$p_j=m_j/m_{j+1}$.
As for $\qbar_j$, using
Theorem~\ref{th:enriques} and the identity (\ref{eq:betabar}), we have
\[
\begin{split}
  \frac{m_{j-1}}{m_{j+1}}\,\qbar_{j-1}-p_j+q_j
  &= \frac{m_{j-1}}{m_{j+1}}\frac{\betabar_{j-1}}{m_j}
  -\frac{m_j}{m_{j+1}}+\frac{\beta_j-\beta_{j-1}+m_j}{m_{j+1}}\\
  &= \frac{1}{m_{j+1}}
  \bigg(
  \frac{m_{j-1}}{m_j}\,\betabar_{j-1}+\beta_j-\beta_{j-1}
  \bigg)\\
  &= \frac{\betabar_j}{m_{j+1}}.
\end{split}
\]
The result follows.
\qed

\section{The proof of Proposition~\ref{p:unibranch}}
\label{s:four}

Let $\xi$ be a jumping number and set $E=\sum_{\beta\in\BBB}E_\beta$
with $\BBB\subset\VVV$ the subset of all vertices such that
$\floor{\xi e_\beta}=\xi e_\beta$.
Then
\begin{equation} \label{eq:inclusion}
  \mu_\ast\Oo_Y(K_\mu-\floor{\xi\mu^\ast C})
  \varsubsetneq  \mu_\ast\Oo_Y(K_\mu-\floor{\xi\mu^\ast C}+E).
\end{equation}
We will show that the right hand side ideal can be computed using only
$\RRR$-chains contained in $E$.  A chain $\Gamma$ of exceptional
divisors $E_\beta$ in the dual graph of the log resolution is called
an $\RRR$-chain if its both extremities are relevant divisors.  A
relevant divisor will be considered as an improper $\RRR$-chain.

\paragraph{Claim}
\begin{equation} \label{eq:equality}
    \mu_\ast\Oo_Y(K_\mu-\floor{\xi\mu^\ast C}+E)
  = \mu_\ast\Oo_Y(K_\mu-\floor{\xi\mu^\ast C}+\sum \Gamma).
\end{equation}
where $\Gamma$ runs through all the maximal $\RRR$-chains contained in
$E$. 

Indeed, suppose that $E$ is a connected subgraph of the dual graph
that contains at least two irreducible components.  Let $E_0$ be an
irreducible component which is not a relevant divisor and which is an
extremity for $E$. %
The intersection of $E_0$ with
$E'=E-E_0$ consists of exactly one point $P$ that will be seen as a
divisor of $E_0$.  Notice that
\begin{equation} \label{eq:degrees}
  \deg (\Ctilde-E)|_{E_0} \leq 1.
\end{equation}
Since $(K_\mu+E)|_{E_0}\sim K_{E_0}+(E-E_0)|_{E_0}=K_{E_0}+P$, we have
the commutative diagram
\[
\begin{CD}
  &&&& 0 && 0\\
  &&&&  @VVV  @VVV \\
  0 @>>> \Oo_Y(\Delta) @>>> \Oo_Y(\Delta+E') @>>>
  \Oo_{E'}(K_{E'}-\floor{\xi\mu^\ast C}|_{E'})@>>> 0 \\
  && @|  @VVV  @VVV \\
  0 @>>> \Oo_Y(\Delta) @>>> \Oo_Y(\Delta+E) @>>>
  \Oo_E(K_E-\floor{\xi\mu^\ast C}|_E)@>>> 0  \\
  &&&&  @VVV  @VVV \\
  &&&& \Oo_{E_0}(\Delta_0+P) @=
  \Oo_{E_0}(\Delta_0+P) \\
  &&&&  @VVV  @VVV \\
  &&&& 0 && 0\\
\end{CD}
\]
where $\Delta=K_\mu-\floor{\xi\mu^\ast C}$ and
$\Delta_0=K_{E_0}-\floor{\xi\mu^\ast C}|_{E_0}$.  The last entry in
the middle vertical short exact sequence is given by the snake lemma.
Pushing down to $X$ this exact sequence and using the local vanishing,
we get that
\[
  \mu_\ast\Oo_Y(K_\mu-\floor{\xi\mu^\ast C}+E)
  = \mu_\ast\Oo_Y(K_\mu-\floor{\xi\mu^\ast C}+E')
\]
if and only if $h^0(E_0,K_{E_0}+P-\floor{\xi\mu^\ast C}|_{E_0})=0$.
But this is true since $E_0$ is an extremity of $E$ and 
\[
\begin{split}
  \floor{\xi\mu^\ast C} \cdot E_0
  &= \sum_{\beta\in\BBB}\xi e_\beta E_\beta\cdot E_0
  + \sum_{\alpha\not\in\BBB} \floor{\xi e_\alpha}E_\alpha\cdot E_0\\
  &> \sum_{\beta\in\BBB}\xi e_\beta E_\beta\cdot E_0
  + \sum_{\alpha\not\in\BBB}(\xi e_\alpha-1)E_\alpha\cdot E_0
  + (\xi-1)\Ctilde\cdot E_0\\
  &= \xi\mu^\ast C \cdot E_0 -
  \big( \Ctilde+\sum_{\alpha\not\in\BBB}E_\alpha \big) \cdot E_0 \\
  &= -(\Ctilde-E) \cdot E_0, 
\end{split}
\]
and with (\ref{eq:degrees}) yield
\[
  \deg P-\floor{\xi\mu^\ast C} \cdot E_0
  < 1 + (\Ctilde-E)\cdot E_0 \leq 2.
\]
Hence
\[
  \deg(K_{E_0}+P-\floor{\xi\mu^\ast C}|_{E_0}) \leq -1.
\]
Repeated use of this argument shows that we can eliminate from $E$ one
component at a time replacing $E$ by $E'$ as long as $E_0$ is not
a relevant divisor.  Hence if $E$ is
connected, only its maximal $\RRR$-chain counts in
computing the ideal $\Oo_Y(K_\mu-\floor{\xi\mu^\ast C}+E)$.  
The above argument, applied to each connected part of $E$, justifies the
claim.  

\medskip

To end the proof of the proposition it is sufficient to show that
the proper $\RRR$-chains can be discarded in the equality
(\ref{eq:equality}).  This is done next.

\paragraph{Claim}
If $\Gamma$ is a proper $\RRR$-chain, then
$h^0(\Gamma,K_\Gamma-\floor{\xi\mu^\ast C}|_\Gamma)=0$.  

Indeed, let us suppose that the $\RRR$-chain $\Gamma$ connects
$E_\rho$ and $E_{\rho'}$, with $\rho'>\rho$.  The first 
$E_\gamma\neq E_\rho$ belonging to $\Gamma$ that $E_\rho$ intersects
is either $E_{\rho+1}$ if $a'_1\geq 2$, or $E_{\rho+a'_2+1}$ if not.
In the former case $e_{\rho+1}=e_\rho+\pi$ and in the latter we have
$e_{\rho+a'_2+1}=(a'_2+1)e_\rho+\pi$.  Here $\pi$ equals a certain
$m_j$, depending on $\rho$.  Since $\xi e_\gamma$ is an integer for
any $E_\gamma\subset\Gamma$, these equalities imply that $\xi\pi$ is
an integer.  Then
\[
  \floor{\xi\mu^\ast C} \cdot E_\rho 
  = \xi\mu^\ast C \cdot E_\rho - \fpart{\xi\mu^\ast C} \cdot E_\rho
  = -\big(\!
  \fpart{\xi e_{\rho-a_m}}E_{\rho-a_m}+\fpart{\xi e_{\rho-1}}E_{\rho-1}
  \big) \cdot E_\rho = 0
\]
as $\pi$ divides both $e_{\rho-a_m}$ and $e_{\rho-1}$ by
Lemma~\ref{l:eCoeffBefore_r}.  Now, if $\Gamma'=\Gamma-E_\rho$ and $P$ is
the intersection point of $E_\rho$ with $\Gamma'$, from the short exact
sequence
\[
  0 \lra
  \Oo_{\Gamma'}(K_{\Gamma'}-\floor{\xi\mu^\ast C}|_{\Gamma'}) \lra
  \Oo_\Gamma(K_\Gamma-\floor{\xi\mu^\ast C}|_\Gamma) \lra
  \Oo_{E_\rho}(K_{E_\rho}+P-\floor{\xi\mu^\ast C}|_{E_\rho}) \lra
  0
\]
we obtain that 
\[
  h^0(\Gamma,K_\Gamma-\floor{\xi\mu^\ast C}|_\Gamma) 
  = h^0(\Gamma',K_{\Gamma'}-\floor{\xi\mu^\ast C}|_{\Gamma'}).
\]
Inductively we cut off from the new chain the lowest extremity
$E_\gamma$ and keep denoting the resulting chain by $\Gamma'$.  We
eventually arrive at
\begin{multline*}
  0 \lra
  \Oo_{E_{\rho'}}(K_{E_{\rho'}}-\floor{\xi\mu^\ast C}|_{E_{\rho'}}) \lra
  \Oo_{\Gamma'}(K_{\Gamma'}-\floor{\xi\mu^\ast C}|_{\Gamma'}) \\ \lra
  \Oo_{E_{\gamma}}(K_{E_{\gamma}}+P-
  \floor{\xi\mu^\ast C}|_{E_{\gamma}}) \lra
  0,
\end{multline*}
where $\gamma$ equals either $\rho'-1$ or $\rho'-a_{m'}$, \ie
$E_\gamma$ is the last exceptional divisor in $\Gamma$ different from
$E_{\rho'}$ and $P$ is the intersection point between $E_\gamma$ and
$E_{\rho'}$.  Arguing as before 
$\deg\floor{\xi\mu^\ast C}|_{E_{\gamma}}=0$.  As for the computation
of $\deg\floor{\xi\mu^\ast C}|_{E_{\rho'}}$, we have
\[
  \floor{\xi\mu^\ast C} \cdot E_{\rho'} 
  = -\Big(\!
  \fpart{\xi e_{\rho'-a'_{m'}}}E_{\rho'-a'_{m'}}
  +\fpart{\xi e_{\rho'-1}}E_{\rho'-1}
  +\fpart{\xi e_{\rho'+\alpha}}E_{\rho'+\alpha}
  \Big) \cdot E_{\rho'}.
\]
Since, as it has been said, $\gamma$ equals either $\rho'-1$ or
$\rho'-a_{m'}$, and since 
$e_{\rho'} = e_{\rho'-a_{m'}}+e_{\rho'-1}+\pi'$ with $\pi'|\pi$, 
\[
  \floor{\xi\mu^\ast C} \cdot E_{\rho'} 
  = \fpart{\xi e_{\rho'+\alpha}} = 0.
\]
Hence $h^0(\Gamma',K_{\Gamma'}-\floor{\xi\mu^\ast C}|_{\Gamma'})=0$
and finally,  
$h^0(\Gamma,K_\Gamma-\floor{\xi\mu^\ast C}|_{\Gamma})=0$
for the proper $\RRR$-chain $\Gamma$, justifying the claim.

\medskip

Now, the short exact sequence 
\begin{multline*}
  0 \lra 
  \mu_\ast\Oo_Y(K_\mu-\floor{\xi\mu^\ast C}) \lra
  \mu_\ast\Oo_Y(K_\mu-\floor{\xi\mu^\ast C}+\sum_{\beta\in\BBB}E_\beta)
  \\ \lra
  \bigoplus_{\rho\in(\BBB\cap\RRR)'}
  H^0(E_\rho,K_{E_\rho}-\floor{\xi\mu^\ast C}|_{E_\rho}) \lra
  0
\end{multline*}
shows that $\xi$ must be contributed by a relevant divisor.  The
set $(\BBB\cap\RRR)'$ is the set of relevant positions in $\BBB$ that
do not define proper $\RRR$-chains contained in $\BBB$.

\bigskip

\begin{flushleft}
  Daniel {\sc Naie}\\
  D\'epartement de Math\'ematiques\\
  Universit\'e d'Angers\\
  F-40045 Angers\\
  France\\
  \small{Daniel.Naie@univ-angers.fr}
\end{flushleft}

\end{document}